# Incomplete Yamabe flows and removable singularities

Mario B. Schulz

*ETH Zürich, Rämistrasse 101, 8092 Zürich, Switzerland*

5 August 2019

We study the Yamabe flow on a Riemannian manifold of dimension $m \geq 3$ minus a closed submanifold of dimension $n$ and prove that there exists an instantaneously complete solution if and only if $n > \frac{m-2}{2}$. In the remaining cases $0 \leq n \leq \frac{m-2}{2}$ including the borderline case, we show that the removability of the $n$-dimensional singularity is necessarily preserved along the Yamabe flow. In particular, the flow must remain geodesically incomplete as long as it exists. This is contrasted with the two-dimensional case, where instantaneously complete solutions always exist.

## Contents

# 1. Introduction

Let $(M, g_0)$ be any Riemannian manifold of dimension $m \geq 3$. We do not necessarily assume that $M$ is compact or complete. However, we always implicitly assume that manifolds and Riemannian metrics are smooth. A family $(g(t))_{t \in [0,T[}$ of Riemannian metrics on $M$ is called *Yamabe flow* with initial metric $g_0$ if

$$\begin{cases} \frac{\partial}{\partial t} g(t) = -\mathrm{R}_{g(t)}\, g(t) & \text{in } M \times [0, T[, \\ g(0) = g_0 & \text{on } M, \end{cases} \tag{1}$$

where $\mathrm{R}_{g(t)}$ denotes the scalar curvature of the Riemannian manifold $(M, g(t))$. Richard Hamilton [9] introduced the Yamabe flow in 1989 as an alternative approach to the Yamabe Theorem which states that any closed Riemannian manifold $(M, g)$ admits a conformal metric of constant scalar curvature. We call a Riemannian metric $\tilde{g}$ on $M$ *conformal* to $g$ if there exists a positive function $u \in \mathrm{C}^\infty(M)$ such that $\tilde{g} = ug$ and we call a manifold *closed* if it is compact without boundary. The Yamabe Theorem was proved by Richard Schoen [20] in 1984 after work of Yamabe [27], Trudinger [26] and Aubin [1] based on the theory of elliptic equations.

Hamilton [9] proved that given any closed Riemannian manifold $(M, g_0)$, there exists a unique solution $(g(t))_{t \in [0,T[}$ to equation (1) for some $T > 0$. On noncompact manifolds however, the well-posedness of the Yamabe flow is more subtle especially in the case that $(M, g_0)$ is geodesically incomplete.

On arbitrary two-dimensional manifolds the Yamabe flow coincides with the Ricci flow and the well-posedness theory of Gregor Giesen and Peter Topping [24, 7, 6, 25] applies. We summarise their results in the following statement.

**Theorem 1** (Giesen–Topping [6], Topping [25])**.** *Let $(M, g_0)$ be any smooth, connected surface. Then there exists a unique solution $(g(t))_{t \in [0,T[}$ of equation (1) on $M$ satisfying*

(I) $g(0) = g_0$,

(II) $(M, g(t))$ *is geodesically complete for every $t > 0$.*

Remarkably, the initial surface $(M, g_0)$ is *not* assumed to be compact or complete. Property (II) is called *instantaneous completeness*. Uniqueness fails without restriction to the class of instantaneously complete flows, as the following example shows.

**Example 1** (Topping [24])**.** Let $g_0$ be the Euclidean metric restricted to

(a) the unit disc $M = \{x \in \mathbb{R}^2 \mid |x| < 1\}$,

(b) the punctured plane $M = \mathbb{R}^2 \setminus \{0\}$.

In both cases (a) and (b), the Riemannian manifold $(M, g_0)$ is geodesically incomplete and $(\bar{g}(t))_{t\in[0,\infty[}$ given by $\bar{g}(t) = g_0$ for all $t \geq 0$ is a constant Yamabe flow on $M$ because $\mathrm{R}_{g_0}$ vanishes identically. However, this flow is not unique. In either case, there exists an instantaneously complete Yamabe flow $(g(t))_{t\in[0,\infty[}$ on $M$ with $g(0) = g_0$ and according to Theorem 1, this flow is unique in the class of Yamabe flows on $M$ satisfying (I) and (II).

## 1.1. Main result

In [21, 22, 23] the author investigated the existence and uniqueness of Yamabe flows on noncompact manifolds of dimension $m \geq 3$ and proved that instantaneously complete Yamabe flows exist unconditionally if the initial manifold is conformal to hyperbolic space. In general however, and in contrast with the two-dimensional case, it is possible that problem (1) does *not* allow any instantaneously complete solution if the initial manifold $(M, g_0)$ is incomplete. A conformally flat example is the punctured sphere: According to [21, Theorem 3], any Yamabe flow starting from the punctured round sphere of dimension $m \geq 3$ must remain geodesically incomplete as long as the flow exists. In this article we generalise this observation to arbitrary manifolds of dimension $m \geq 3$ and show that the statement is still true if we remove not only a point but a low-dimensional submanifold. More precisely, we prove the following sharp result.

**Theorem 2.** *Let $(M, g_0)$ be a closed Riemannian manifold of dimension $m \geq 3$ and let $\emptyset \neq N \subset M$ be a closed submanifold of dimension $n \geq 0$.*

(i) *If $n > \frac{m-2}{2}$ then an instantaneously complete Yamabe flow $(g(t))_{t\in[0,\infty[}$ on $M \setminus N$ with $g(0) = g_0$ exists.*

(ii) *If $n \leq \frac{m-2}{2}$ then any Yamabe flow with initial data $(M \setminus N, g_0)$ is incomplete and uniquely given by the restriction of the Yamabe flow with initial data $(M, g_0)$.*

*Remark.* Theorem 2 states that the removability of $n$-dimensional singularities is necessarily preserved along the Yamabe flow on closed manifolds of dimension $m \geq 3$ if and only if $n \leq \frac{m-2}{2}$. The borderline case $n = \frac{m-2}{2}$ for even $m \geq 4$ is particularly delicate. Moreover, Theorem 2 does *not* extend to the case $m = 2$ and $n = 0$ because on surfaces instantaneously complete Yamabe flows always exist according to Theorem 1. The singular submanifold $N$ is not assumed to be connected: If $n = 0$ then $N \subset M$ is a finite set of points. In case (i) the uniqueness in the class of instantaneously complete Yamabe flows with initial metric $g_0$ is still open in general.

The Yamabe flow has been studied on incomplete manifolds by Bahuaud and Vertman [4, 3] who proved short-time existence and uniqueness of solutions to equation (1) in a suitable class of flows which preserve the incompleteness. In Theorem 2 (ii), we allow any flow satisfying the initial condition and do not require further restrictions.

## 1.2. Similarities with the singular Yamabe problem

Theorem 2 can be interpreted as parabolic analogue of several known result about elliptic (Yamabe-type) equations. Most notably, the threshold $\frac{m-2}{2}$ also appears in the "singular Yamabe Problem" which we formulate on the sphere:

**Theorem 3** (Singular Yamabe problem [12, 19, 18, 14, 15]). *Let $(\mathbb{S}^m, g_{\mathbb{S}^m})$ be the standard sphere of dimension $m \geq 3$ and let $N \subset \mathbb{S}^m$ be a closed submanifold of dimension $n \geq 1$.*

(i) *If $n > \frac{m-2}{2}$ then there exists a complete, conformal metric on $\mathbb{S}^m \setminus N$ with constant negative scalar curvature.*

(ii) *If $n \leq \frac{m-2}{2}$ then there exists a complete, conformal metric on $\mathbb{S}^m \setminus N$ with constant positive scalar curvature.*

Loewner and Nirenberg [12] studied case (i) of Theorem 3. Moreover, they proved [12, Theorem 7] that given any open domain $\Omega \subset \mathbb{R}^m$, compact subsets $N \subset \Omega$ of Hausdorff dimension $n < \frac{m-2}{2}$ are removable sets for positive solutions of the equation

$$\Delta u = u^{\frac{m+2}{m-2}}$$

and conjectured that the result extends to the borderline case $n = \frac{m-2}{2}$. For the proof of Theorem 2 it is relevant that the statement of Theorem 3 (i) holds on any compact Riemannian manifold in place of the sphere:

**Theorem 4** (Aviles–McOwen [2]). *Let $(M, g)$ be a compact Riemannian manifold of dimension $m \geq 3$ and let $N \subset M$ be a closed submanifold of dimension $n$. Then there exists a complete, conformal metric $\hat{g}$ on $M \setminus N$ with constant negative scalar curvature if and only if $n > \frac{m-2}{2}$.*

On the other hand, Schoen and Yau [18, Theorem 2.7] showed that if $M \subset \mathbb{S}^m$ is a (dense) domain in the sphere which has a complete, conformal metric $g$ with scalar curvature $\mathrm{R}_g \geq 1$, then the Hausdorff dimension of the boundary $\partial M$ is at most $\frac{m-2}{2}$. Hence, the condition $n \leq \frac{m-2}{2}$ is necessary in Theorem 3 (ii). Moreover, Schoen [19] studied the case where $N \subset \mathbb{S}^m$ is a finite set of at least two points. In the setting of Theorem 3 (ii), Pacard [16] considered the borderline case $n = \frac{m-2}{2}$ for even $m$ first and then solved the singular Yamabe Problem for $n \leq \frac{m-2}{2}$ in joint work with Mazzeo [14, 15].

A variant of the singular Yamabe problem asks for complete conformal metrics with vanishing scalar curvature. Ma and McOwen [13] generalised work of Lee and Parker [11] and Jin [10] and proved that if $(M, g)$ is a compact Riemannian manifold of dimension $m \geq 3$ with positive Yamabe invariant and if $N \subset N$ is a closed submanifold of dimension $n \leq \frac{m-2}{2}$ then $M \setminus N$ admits a complete conformal metric with constant zero scalar curvature.

## 1.3. Geometric lemmata

The conformal class of the initial metric $g_0$ is preserved under the Yamabe flow. In fact, any solution $(g(t))_{t\in[0,T[}$ of problem (1) is of the form $g(t) = u(\cdot,t)g_0$ with conformal factor $u\colon M \times [0,T[ \to ]0,\infty[$ satisfying

$$\begin{cases} \frac{\partial}{\partial t}u = -\mathrm{R}_g u & \text{in } M \times ]0,T[, \\ u = 1 & \text{on } M \times \{0\}. \end{cases} \qquad (2)$$

**Lemma 1.1.** *Let $(M,g)$ be a Riemannian manifold of dimension $m \geq 3$ and let $\tilde{g} = U^{\frac{4}{m-2}}g$ be a conformal metric on $M$. Then the scalar curvature of $(M,\tilde{g})$ is*

$$\mathrm{R}_{\tilde{g}} = U^{-\frac{m+2}{m-2}}\Big(\mathrm{R}_g U - 4\tfrac{m-1}{m-2}\Delta_g U\Big).$$

One of many references containing a proof of Lemma 1.1 is [17, §18]. Depending on the dimension $m \geq 3$ of the manifold $M$ we introduce the exponent

$$\eta := \frac{m-2}{4} \qquad (3)$$

and apply Lemma 1.1 in (2) to show that $U = u^\eta$ solves the parabolic equation

$$\frac{1}{m-1}\frac{\partial U}{\partial t} = \Big(-\frac{m-2}{4(m-1)}\mathrm{R}_{g_0}U + \Delta_{g_0}U\Big)U^{-\frac{1}{\eta}} \qquad (4)$$

or equivalently,

$$\frac{1}{\eta+1}\frac{\partial U^{1+\frac{1}{\eta}}}{\partial t} = -\mathrm{R}_{g_0}U + \frac{m-1}{\eta}\Delta_{g_0}U. \qquad (5)$$

**Lemma 1.2.** *Let $(M,g)$ be a Riemannian manifold of dimension $m$ and let $N \subset M$ be a closed submanifold of dimension $0 \leq n < m$. Let $r\colon M \to [0,\infty[$ measure the Riemannian distance from $N$. Then, as $r \to 0$,*

$$r\Delta_g r = m - n - 1 + O(r).$$

*Proof.* The convergence $r\Delta_g r \to m - n - 1$ as $r \to 0$ is shown in [12, (5.5)] for the case that $(M,g_0)$ is Euclidean space and can also be verified in this more general setting by a similar computation. Since $r^2$ is smooth in some neighbourhood of $N$ and $\Delta_g r^2 = 2 + 2r\Delta_g r$ we may conclude that the error term $O(r)$ is bounded uniformly by a function linear in $r$. □

The next lemma is related to results by Escobar [5] but we provide a self-contained proof. In the following, we say that two Riemannian metrics $\bar{g}$ and $g_0$ on $M$ are *comparable* if there exist constants $c_2 \geq c_1 > 0$ such that $c_1\bar{g} \leq g_0 \leq c_2\bar{g}$.

**Lemma 1.3.** *Let $(M, g_0)$ be any Riemannian manifold of dimension $m \geq 3$ and let $N \subset M$ be a closed submanifold of dimension $0 \leq n < m$. Then there exist a neighbourhood $\Omega \subset M$ of $N$ and a Riemannian metric $\overline{g}$ on $M$ which is conformal and comparable to $g_0$ such the scalar curvature of $(M, \overline{g})$ vanishes in $\Omega$.*

*Proof.* Let $B_\varepsilon(x)$ be the metric ball of radius $0 < \varepsilon \leq 1$ around $x$ in $(M, g_0)$ and let

$$\Omega_\varepsilon := \bigcup_{x \in N} B_\varepsilon(x).$$

According to Lemma 1.1 the scalar curvature of $g = U^{\frac{4}{m-2}} g_0$ vanishes in $\Omega_\varepsilon$ if

$$\begin{cases} -\Delta_{g_0} U + \frac{m-2}{4(m-1)} \mathrm{R}_{g_0} U = 0 & \text{in } \Omega_\varepsilon, \\ \phantom{-\Delta_{g_0} U + \frac{m-2}{4(m-1)} \mathrm{R}_{g_0}} U = 1 & \text{on } \partial \Omega_\varepsilon. \end{cases} \qquad (6)$$

We consider the affine space $X = 1 + H_0^1(\Omega_\varepsilon)$ and the energy $E \colon X \to \mathbb{R}$ given by

$$E(U) = \int_{\Omega_\varepsilon} |\nabla U|_{g_0}^2 + \tfrac{m-2}{4(m-1)} \mathrm{R}_{g_0} U^2 \, d\mu_{g_0}.$$

By Young's inequality, there exists a finite constant $C > 0$ depending only on $\inf_{\Omega_1} \mathrm{R}_{g_0}$ and $m$ but not on $\varepsilon$ or $U$ such that

$$E(U) \geq \int_{\Omega_\varepsilon} |\nabla U|_{g_0}^2 \, d\mu_{g_0} - C(U-1)^2 - C \, d\mu_{g_0}.$$

We claim that $E$ is coercive if $\varepsilon > 0$ is sufficiently small. Let

$$\lambda_\varepsilon = \inf_{u \in H_0^1(\Omega_\varepsilon)} \frac{\int_{\Omega_\varepsilon} |\nabla u|_{g_0}^2 \, d\mu_{g_0}}{\int_{\Omega_\varepsilon} u^2 \, d\mu_{g_0}}$$

be the smallest Dirichlet eigenvalue of $-\Delta_{g_0}$ on $\Omega_\varepsilon$. Given any $U \in X$ we have $U - 1 \in H_0^1(\Omega_\varepsilon)$ and hence

$$E(U) \geq \int_{\Omega_\varepsilon} \left(1 - \frac{C}{\lambda_\varepsilon}\right) |\nabla U|_{g_0}^2 - C \, d\mu_{g_0}.$$

Coercivity of $E$ for sufficiently small $\varepsilon > 0$ follows provided that $\lambda_\varepsilon \to \infty$ as $\varepsilon \searrow 0$. Indeed, if $u_\varepsilon \in H_0^1(\Omega_\varepsilon)$ is a corresponding eigenfunction which is normalised to $\|u_\varepsilon\|_{L^2(\Omega_\varepsilon)}^2 = 1$ such that $\lambda_\varepsilon = \|\nabla u_\varepsilon\|_{L^2(\Omega_\varepsilon)}^2$ then by Hölder's and Sobolev's inequality

$$1 = \int_{\Omega_\varepsilon} u_\varepsilon^2 \, d\mu_{g_0} \leq \left(\int_{\Omega_\varepsilon} d\mu_{g_0}\right)^{\frac{2}{m}} \left(\int_{\Omega_\varepsilon} |u_\varepsilon|^{\frac{2m}{m-2}} \, d\mu_{g_0}\right)^{\frac{m-2}{m}} \leq \left(\int_{\Omega_\varepsilon} d\mu_{g_0}\right)^{\frac{2}{m}} S \lambda_\varepsilon. \qquad (7)$$

The Sobolev constant $S$ does not depend on $0 < \varepsilon \leq \varepsilon_0$ if we extend each $u_\varepsilon$ by zero to $\Omega_{\varepsilon_0}$ and apply Sobolev's inequality there. The volume of $\Omega_\varepsilon$ however converges to zero as $\varepsilon \searrow 0$ and $\lambda_\varepsilon \to \infty$ follows according to (7).

By the argument above, there exists some $\varepsilon > 0$ such that $E$ is coercive. Since $E$ is also weakly sequentially lower semi-continuous, $E$ attains its minimum at some $U \in X$ which satisfies the corresponding Euler-Lagrange equation (6). Since $E(U) = E(|U|)$ we may assume $U \geq 0$. By elliptic regularity theory, $U$ is smooth up to the boundary and strictly positive by the strong maximum principle. To prove the claim, we extend $U$ to a uniformly bounded, positive function on all of $M$ and set

$$\overline{g} = U^{\frac{4}{m-2}} g.$$
□

## 2. Boundedness

Theorem 2 (ii) is based on the following local upper bound for singular Yamabe flows.

**Theorem 5.** *Let $(M, g_0)$ be any open Riemannian manifold of dimension $m \geq 3$ and let $N \subset M$ be any closed submanifold of dimension $0 \leq n \leq \frac{m-2}{2}$. Let $(g(t))_{t \in [0,T]}$ be any Yamabe flow on $M \setminus N$ with $g(0) = g_0$. Then there exists a constant $C > 0$ such that $g(t) \leq C g_0$ holds for all $t \in [0, T]$ in some neighbourhood of $N$.*

*Remark.* Theorem 5 already implies that every Yamabe flow $(g(t))_{t \in [0,T]}$ on $M \setminus N$ with $g(0) = g_0$ must be geodesically incomplete for all $t \in [0, T]$. The intuition for Theorem 5 is that conformally opening up the edge $N$ generates positive scalar curvature near $N$ which prevents the Yamabe flow from expanding further. In fact, our argument is based on the construction of such a conformal background metric.

Since we allow the manifold $M$ to be noncompact and do not make any further assumption, it is unknown whether any Yamabe flow on $M$ with initial metric $g_0$ exists. Therefore, existence is not part of the claim of Theorem 5.

### 2.1. The key idea

Let $(M, g_0)$ be an open Riemannian manifold of dimension $m \leq 3$ and let $N \subset M$ be a closed submanifold of dimension $n \leq \frac{m-2}{2}$. According to Lemma 1.3, there exist a neighbourhood $\Omega \subset M$ of $N$ and a Riemannian metric $\overline{g}$ on $M$ which is conformal and comparable to $g_0$ such the scalar curvature of $(M, \overline{g})$ vanishes in $\Omega$. Let $B_\delta(x)$ denote the metric ball of radius $\delta > 0$ around the point $x$ in $(M, \overline{g})$ and let

$$\Omega_\delta := \bigcup_{x \in N} B_\delta(x).$$

Let $r \colon M \to [0, \infty[$ be the distance from $N$ in $(M, \overline{g})$. Since $N$ is smooth and closed, there exists some $\delta > 0$ such that $r^2$ restricted to $\Omega_\delta$ is smooth. We may assume that $\delta > 0$ is sufficiently small such that $\Omega_\delta \subset \Omega$ which means that $\mathrm{R}_{\overline{g}} = 0$ in $\Omega_\delta$.

Let $f \in C^\infty(M \setminus N)$ be a positive function to be chosen, let $\eta = \frac{m-2}{4}$ as in (3) and let

$$\tilde{g} = f\bar{g} \quad \text{on } M \setminus N. \tag{8}$$

Since $R_{\bar{g}} = 0$ in $\Omega_\delta$, Lemma 1.1 applied to (8) respectively $\bar{g} = f^{-1}\tilde{g}$ yields

$$-f^{-\eta-1}\Delta_{\bar{g}}f^\eta = \frac{\eta R_{\tilde{g}}}{m-1} = f^\eta \Delta_{\tilde{g}}f^{-\eta} \quad \text{in } \Omega_\delta \setminus N. \tag{9}$$

Let $(g(t))_{t \in [0,T]}$ be any Yamabe flow on $M \setminus N$ with $g(0) = g_0$. Since $\tilde{g}$ is conformal to $g_0$ on $M \setminus N$, there exists a function $U\colon (M \setminus N) \times [0,T] \to {]0,\infty[}$ such that

$$g(t) = U(\cdot, t)^{\frac{1}{\eta}}\tilde{g}.$$

We combine (4) and (9) to obtain

$$\begin{aligned}
\frac{1}{m-1}\frac{\partial U}{\partial t} &= \left(-\frac{\eta R_{\tilde{g}}}{m-1}U + \Delta_{\tilde{g}}U\right)U^{-\frac{1}{\eta}} \\
&= \left(-\frac{\eta R_{\tilde{g}}}{m-1}(U - cf^{-\eta}) + \Delta_{\tilde{g}}(U - cf^{-\eta})\right)U^{-\frac{1}{\eta}}
\end{aligned} \tag{10}$$

in $(\Omega_\delta \setminus N) \times [0,T]$ with an arbitrary constant $c \in \mathbb{R}$. Equation (10) is the central step in the proof and reveals the key idea: If we show $U(\cdot, t) \leq cf^{-\eta}$ for some $c > 0$ and all $t \in [0,T]$, then, since $\bar{g}$ and $g_0$ are comparable, we obtain

$$g(t) \leq c^{\frac{1}{\eta}}f^{-1}\tilde{g} = c^{\frac{1}{\eta}}\bar{g} \leq Cg_0 \quad \text{in } \Omega_\delta \setminus N. \tag{11}$$

## 2.2. Low-dimensional singularities

To implement the key idea, we require the metric $\tilde{g}$ defined in (8) to be geodesically complete towards $N$ with positive scalar curvature in $\Omega_\delta \setminus N$. This is easier to achieve if the dimension $n$ of the submanifold $N$ is strictly below the threshold $\frac{m-2}{2}$. Choosing $f = r^{-2}$ in $\Omega_\delta \setminus N$, equation (9) implies

$$R_{\tilde{g}} = -\frac{m-1}{\eta}r^{2\eta+2}\Delta_{\bar{g}}r^{-2\eta} = -(m-1)\big(m - 2r\Delta_{\bar{g}}r\big) \tag{12}$$

in $\Omega_\delta \setminus N$, where we computed $2(2\eta + 1) = m$. By Lemma 1.2,

$$R_{\tilde{g}} = (m-1)\big(m - 2 - 2n + 2O(r)\big) \tag{13}$$

as $r \to 0$. If $\delta > 0$ is sufficiently small and $n < \frac{m-2}{2}$ then (13) is positive in $\Omega_\delta \setminus N$.

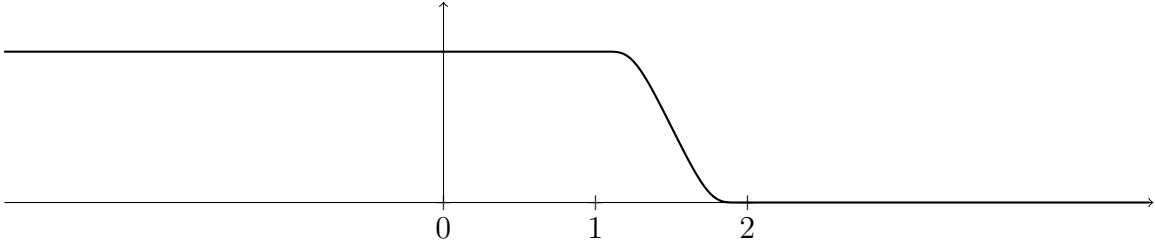

Figure 1: Graph of a suitable cutoff function $\chi$.

Let $\chi\colon \mathbb{R} \to [0,1]$ be a smooth, nonincreasing cutoff function satisfying $\chi(s) = 1$ for $s \leq 1$ and $\chi(s) = 0$ for $s \geq 2$ as shown in Figure 1. We recall $r = \operatorname{dist}_{\overline{g}}(x, N)$ and $\tilde{g} = f\,\overline{g}$ with $f = r^{-2}$. Given $\varepsilon > 0$ the function $\rho\colon \Omega_\delta \setminus N \to \,]0, \infty[$ defined by

$$\rho(x) = \varepsilon \operatorname{dist}_{\tilde{g}}(x, \partial\Omega_\delta) \tag{14}$$

satisfies $\frac{\partial}{\partial r}\rho = -\varepsilon\sqrt{f}$ and $\rho(x) \to \infty$ as $r(x) \to 0$. Let $\varphi\colon \Omega_a \to [0,1]$ be defined by

$$\varphi = \chi \circ \rho. \tag{15}$$

Using the formulas from Lemma A.4 in the appendix, we have

$$\begin{aligned}
\frac{\partial \varphi}{\partial r} &= -\varepsilon\sqrt{f}\,\chi' \circ \rho, \\
\frac{\partial^2 \varphi}{\partial r^2} &= \varepsilon^2 f\,\chi'' \circ \rho - \varepsilon \frac{\partial \sqrt{f}}{\partial r} \chi' \circ \rho, \\
|\nabla \varphi|^2_{\tilde{g}} &= \frac{1}{f}\left|\frac{\partial \varphi}{\partial r}\right|^2 = \varepsilon^2 |\chi' \circ \rho|^2, \\
\Delta_{\tilde{g}} \varphi &= \frac{1}{f} \Delta_{\overline{g}} \varphi + \frac{m-2}{2f^2}\langle \nabla f, \nabla \varphi \rangle_{\overline{g}} \\
&= \frac{1}{f}\frac{\partial^2 \varphi}{\partial r^2} + \frac{1}{f}\frac{\partial \varphi}{\partial r}\Delta_{\overline{g}} r + \frac{(m-2)}{2f^2}\frac{\partial f}{\partial r}\frac{\partial \varphi}{\partial r} \\
&= \varepsilon^2 \chi'' \circ \rho - \left((m-1)\frac{1}{f}\frac{\partial \sqrt{f}}{\partial r} + \frac{\Delta_{\overline{g}} r}{\sqrt{f}}\right)\varepsilon \chi' \circ \rho
\end{aligned} \tag{16,17}$$

Since $f = r^{-2}$ the factor

$$-\left((m-1)\frac{1}{f}\frac{\partial \sqrt{f}}{\partial r} + \frac{\Delta_{\overline{g}} r}{\sqrt{f}}\right) = (m-1) - r\Delta_{\overline{g}} r = n + O(r) \tag{18}$$

is bounded uniformly in $\Omega_\delta \setminus N$ by Lemma 1.2. In particular, by choosing the cutoff function $\chi$ appropriately (see [21, Lemma A.4]), there exists a constant $C$ such that

$$\left(\frac{2|\nabla\varphi|^2_{\tilde{g}}}{\varphi} - \Delta_{\tilde{g}} \varphi\right) \leq C\varepsilon \varphi^{1-\frac{1}{\eta}}. \tag{19}$$

After fixing $\delta > 0$ such that $\mathrm{R}_{\tilde{g}} > 0$ in $\Omega_\delta \setminus N$, we choose $c > 0$ such that
$$U \leq V := cf^{-\eta} \quad \text{on } (\Omega_\delta \setminus N \times \{0\}) \cup (\partial\Omega_\delta \times [0, T]). \tag{20}$$
This is possible, since the construction
$$g_0 = U(\cdot, 0)^{\frac{1}{\eta}} \tilde{g} = \left(U(\cdot, 0)f^\eta\right)^{\frac{1}{\eta}} \overline{g}$$
implies that $U(\cdot, 0)f^\eta$ is uniformly bounded in $\Omega_\delta \setminus N$. Moreover, we could prove local upper bounds for $U$ as in [22, Lemma 7] to estimate $c$ in terms of the dimension $m$, the final time $T$, the initial data $U(\cdot, 0)$ and $\delta$. The latter two depend only on $g_0$. Based on (10), we compute
$$\begin{aligned}
&\frac{1}{m-1} \frac{\partial (U - V)\varphi}{\partial t} \\
&= \left(-\frac{\eta \mathrm{R}_{\tilde{g}} \varphi}{m-1}(U - V) + \varphi \Delta_{\tilde{g}}(U - V)\right) U^{-\frac{1}{\eta}} \\
&= -\frac{\eta \mathrm{R}_{\tilde{g}} \varphi}{m-1}(U - V) U^{-\frac{1}{\eta}} + U^{-\frac{1}{\eta}} \Delta_{\tilde{g}}\left((U - V)\varphi\right) \\
&\quad - U^{-\frac{1}{\eta}} \left\langle \nabla\left((U - V)\varphi\right), \nabla \varphi \right\rangle_{\tilde{g}} + \left(\frac{2|\nabla \varphi|_{\tilde{g}}^2}{\varphi} - \Delta_{\tilde{g}} \varphi\right) U^{-\frac{1}{\eta}}(U - V).
\end{aligned} \tag{21}$$
The final argument is the same as for [21, Proposition 2.1] but we repeat it for clarity: By construction, the function $w \colon [0, T] \to [0, \infty[$ given by
$$w(t) = \max_{\Omega_\delta \setminus N} \left(\left((U - V)\varphi\right)(\cdot, t)\right)$$
is well defined and for each $t_0 \in \,]0, T]$ there exists an interior point $q_0 \in \Omega_\delta \setminus N$ such that $w(t_0) = \left((U - V)\varphi\right)(q_0, t_0)$. If $w(t_0) > 0$, then
$$\left(U^{-\frac{1}{\eta}}(U - V)\right)(q_0, t_0) \leq \left((U - V)^{1-\frac{1}{\eta}}\right)(q_0, t_0). \tag{22}$$
Only the last term in equation (21) can be positive at $(q_0, t_0)$ since $q_0 \in \Omega_\delta \setminus N$ is an interior maximum of $(U - V)\varphi(\cdot, t_0)$ and $\mathrm{R}_{\tilde{g}} > 0$ by (13). In combination with estimates (19) and (22) we obtain
$$\begin{aligned}
&\limsup_{\tau \searrow 0} \frac{w(t_0) - w(t_0 - \tau)}{\tau} \\
&\leq \limsup_{\tau \searrow 0} \frac{((U - V)\varphi)(q_0, t_0) - ((U - V)\varphi)(q_0, t_0 - \tau)}{\tau} = \frac{\partial (U - V)\varphi}{\partial t}(q_0, t_0) \\
&\leq (m-1)C\varepsilon \left((U - V)\varphi\right)^{1-\frac{1}{\eta}}(q_0, t_0) = (m-1)C\varepsilon \left(w(t_0)\right)^{1-\frac{1}{\eta}}.
\end{aligned} \tag{23}$$
Since $w(0) = 0$ by (20), Lemma A.2 given in the appendix yields $w(t) \leq \left(\frac{m-1}{\eta}\varepsilon Ct\right)^\eta$ for every $t \in [0, T]$. Letting $\varepsilon \searrow 0$ such that $\varphi \to 1$ pointwise proves $U \leq V$ in $(\Omega_\delta \setminus N) \times [0, T]$ and (11) follows.

## 2.3. The borderline case

In the borderline case $n = \frac{m-2}{2}$ the scalar curvature (13) is not necessarily positive. To overcome this obstruction, we choose the conformal factor $f$ more carefully:

$$f = r^{-2}(-\log r)^{-\frac{2}{3}} \quad \text{in } \Omega_\delta \setminus N.$$

By (9) we have

$$\begin{aligned}
\frac{-R_{\tilde{g}}}{m-1} &= \frac{1}{\eta} f^{-\eta-1} \Delta_{\overline{g}} f^\eta \\
&= f^{-2} \Delta_{\overline{g}} f + (\eta-1) f^{-3} |\nabla f|^2_{\overline{g}} \\
&= f^{-2} \frac{\partial^2 f}{\partial r^2} + f^{-2} \frac{\partial f}{\partial r} \Delta_{\overline{g}} r + \frac{(m-6)}{4} f^{-3} \left(\frac{\partial f}{\partial r}\right)^2.
\end{aligned} \quad (24)$$

Therefore, we compute

$$\frac{\partial f}{\partial r} = \frac{2(1+3\log r)}{3r^3(-\log r)^{\frac{5}{3}}} = \left(\tfrac{2}{3}(-\log r)^{-\frac{1}{3}} - 2(-\log r)^{\frac{2}{3}}\right) rf^2, \quad (25)$$

$$\frac{\partial^2 f}{\partial r^2} = \left(\tfrac{10}{9}(-\log r)^{-\frac{4}{3}} - \tfrac{10}{3}(-\log r)^{-\frac{1}{3}} + 6(-\log r)^{\frac{2}{3}}\right) f^2, \quad (26)$$

$$\begin{aligned}
\left(\frac{\partial f}{\partial r}\right)^2 &= \tfrac{4}{9}(1+3\log r)^2 (-\log r)^{-\frac{4}{3}} f^3 \\
&= 4\left(\tfrac{1}{9}(-\log r)^{-\frac{4}{3}} - \tfrac{2}{3}(-\log r)^{-\frac{1}{3}} + (-\log r)^{\frac{2}{3}}\right) f^3
\end{aligned} \quad (27)$$

and insert (25), (26), (27) in (24) which yields

$$\begin{aligned}
\frac{-R_{\tilde{g}}}{m-1} &= \left(\frac{10}{9} + \frac{(m-6)}{9}\right)(-\log r)^{-\frac{4}{3}} \\
&\quad + \left(-\frac{10}{3} + \frac{2}{3}(r\Delta_{\overline{g}} r) - \frac{2(m-6)}{3}\right)(-\log r)^{-\frac{1}{3}} \\
&\quad + \left(6 - 2(r\Delta_{\overline{g}} r) + (m-6)\right)(-\log r)^{\frac{2}{3}} \\
&= \frac{m+4}{9(-\log r)^{\frac{4}{3}}} - \frac{2(n-O(r))}{3(-\log r)^{\frac{1}{3}}} - \left(m-2-2n+2O(r)\right)(-\log r)^{\frac{2}{3}}
\end{aligned} \quad (28)$$

using the formula $r\Delta_{\overline{g}} r = m - n - 1 + O(r)$ from Lemma 1.2. In the case $n = \frac{m-2}{2}$, the second term in (28) dominates the first and third term as $r \to 0$ which implies that $R_{\tilde{g}} > 0$ in $\Omega_\delta \setminus N$ if $\delta > 0$ is chosen sufficiently small.

Given $0 < \varepsilon < 1$ we define $\rho(x) = \varepsilon \operatorname{dist}_{\tilde{g}}(x, \partial\Omega_\delta)$ and $\varphi = \chi \circ \rho$ as in (14) and (15).

For sufficiently small $\delta > 0$, we may replace (18) by the estimate

$$-\left((m-1)\frac{1}{f}\frac{\partial\sqrt{f}}{\partial r} + \frac{\Delta_{\bar{g}}r}{\sqrt{f}}\right)$$
$$= (m-1)\left((-\log r)^{\frac{1}{3}} - \tfrac{1}{3}(-\log r)^{-\frac{2}{3}}\right) - \left(m-n-1+O(r)\right)(-\log r)^{\frac{1}{3}}$$
$$= \frac{1-m}{3(-\log r)^{\frac{2}{3}}} + \left(n - O(r)\right)(-\log r)^{\frac{1}{3}} \leq n(-\log r)^{\frac{1}{3}}. \tag{29}$$

We now fix $\delta > 0$ such that all the estimates described above are satisfied in $\Omega_\delta \setminus N$. In the support of $\varphi = \chi \circ \rho$ we have by definition

$$2 \geq \rho = \varepsilon \int_r^\delta \frac{1}{s(-\log s)^{\frac{1}{3}}}\, ds = \frac{3\varepsilon}{2}\left((-\log r)^{\frac{2}{3}} - (-\log \delta)^{\frac{2}{3}}\right)$$

and hence $(-\log r)^{\frac{2}{3}} \leq \frac{4}{3\varepsilon} + (-\log \delta)^{\frac{2}{3}}$. Assuming $0 < \varepsilon \leq \frac{2}{3}(-\log \delta)^{-\frac{2}{3}}$ we obtain

$$(-\log r)^{\frac{1}{3}} \leq \sqrt{\frac{2}{\varepsilon}}. \tag{30}$$

We combine (30) and (29) in (17) and recall $\chi' \leq 0$ to conclude

$$\Delta_{\bar{g}}\varphi \geq \varepsilon^2 \chi'' \circ \rho + n\sqrt{2\varepsilon}\, \chi' \circ \rho.$$

Consequently, estimate (19) with $\sqrt{\varepsilon}$ instead of $\varepsilon$ holds also in the borderline case and we may complete the proof of Theorem 5 as before.

## 3. Preservation of removability

In this section, we prove Theorem 2. Our approach for part (ii) is inspired by the energy method Topping [25] used to prove the uniqueness part in Theorem 1 and by the generalisation of this approach to higher dimensions established in [21, § 2.2].

*Proof of Theorem 2.* (i) Let $(M, g_0)$ be a closed Riemannian manifold of dimension $m \geq 3$ and let $N \subset M$ be a closed $n$-dimensional submanifold. By Theorem 4, the condition $n > \frac{m-2}{2}$ is equivalent to the existence of a complete conformal metric $\hat{g}$ on $M \setminus N$ with constant negative scalar curvature. In [22] the author proves the existence of instantaneously complete Yamabe flows on hyperbolic space starting from any given conformally hyperbolic initial metric. The approach is based on an exhaustion of the complete background manifold with smooth, bounded domains and exploits only the negativity of the background scalar curvature and general local upper and lower bounds for Yamabe flows [22, Lemmata 5 and 7]. Therefore, we may replace hyperbolic space by $(M \setminus N, \hat{g})$ and use exactly the same arguments as in [22] to construct an instantaneously complete Yamabe on $M \setminus N$ with initial metric $g_0$.

(ii) Let $(g(t))_{t\in[0,T_0]}$ be a Yamabe flow on $M \setminus N$ with $g(0) = g_0$. Since $M$ is closed, there exists $\tilde{T} > 0$ and a unique Yamabe flow $(\tilde{g}(t))_{t\in[0,\tilde{T}]}$ on $M$ with $\tilde{g}(0) = g_0$ by the results of Hamilton [9]. Let $T = \min\{T_0, \tilde{T}\}$ and let $u, \tilde{u}\colon (M \setminus N) \times [0,T] \to \,]0,\infty[$ be such that $g(t) = u(\cdot,t)g_0$ and $\tilde{g}(t) = \tilde{u}(\cdot,t)g_0$ in $M \setminus N$ for every $t \in [0,T]$.

Since $M$ and $N$ are closed, Theorem 5 implies that $u$ is uniformly bounded from above in the noncompact manifold $(M \setminus N) \times [0,T]$. The same holds for $\tilde{u}$, since it is the restriction of a smooth function on a compact manifold. Additionally, $\inf_{M\times[0,T]} \tilde{u} > 0$. Both, $u^\eta$ and $\tilde{u}^\eta$ are solutions to the equation (5), i.e.

$$\frac{1}{\eta+1}\frac{\partial u^{\eta+1}}{\partial t} = -\mathrm{R}_{g_0}u^\eta + \frac{m-1}{\eta}\Delta_{g_0}u^\eta$$

where $\eta = \frac{m-2}{4}$ as in (3). Let $w := u^{\eta+1} - \tilde{u}^{\eta+1}$ and $w_+(x,t) := \max\{w(x,t), 0\}$. Given $0 \leq \varphi \in \mathrm{C}_c^\infty(M \setminus N)$ we study the evolution of

$$J(t) = \int_M w_+(\cdot,t)\,\varphi\,d\mu_{g_0}.$$

For every fixed $t \in \,]0,T]$ and every $0 < \tau < t$ we have

$$\begin{aligned}
J(t) - J(t-\tau) &= \int_M w_+(\cdot,t)\,\varphi\,d\mu_{g_0} - \int_M w_+(\cdot,t-\tau)\,\varphi\,d\mu_{g_0} \\
&\leq \int_{\{w(\cdot,t)>0\}} \bigl(w_+(\cdot,t) - w_+(\cdot,t-\tau)\bigr)\varphi\,d\mu_{g_0} \\
&\leq \int_{\{w(\cdot,t)>0\}} \bigl(w(\cdot,t) - w(\cdot,t-\tau)\bigr)\varphi\,d\mu_{g_0}
\end{aligned}$$

and obtain

$$\begin{aligned}
\Psi(t) := \limsup_{\tau \searrow 0} \frac{J(t) - J(t-\tau)}{\tau} &\\
\leq \limsup_{\tau \searrow 0} \int_{\{w(\cdot,t)>0\}} \frac{1}{\tau}\bigl(w(\cdot,t) - w(\cdot,t-\tau)\bigr)\varphi\,d\mu_{g_0} &\\
= \int_{\{w(\cdot,t)>0\}} \frac{\partial w}{\partial t}(\cdot,t)\,\varphi\,d\mu_{g_0} &\\
= (\eta+1)\int_{\{w(\cdot,t)>0\}} (u^\eta - \tilde{u}^\eta)(\cdot,t)(-\mathrm{R}_{g_0})\varphi\,d\mu_{g_0} & \quad (31)\\
+ (\eta+1)\frac{m-1}{\eta}\int_{\{w(\cdot,t)>0\}} \varphi\Delta_{g_0}(u^\eta - \tilde{u}^\eta)(\cdot,t)\,d\mu_{g_0}.
\end{aligned}$$

Let $f := (u^\eta - \tilde{u}^\eta)(\cdot,t)$. For any $0 \leq \varphi \in \mathrm{C}_c^\infty(M \setminus N)$ there holds

$$\int_{\{f>0\}} (\varphi\Delta_{g_0}f - f\Delta_{g_0}\varphi)\,d\mu_{g_0} \leq 0 \qquad (32)$$

as shown in [21, § 2.2]. For convenience of the reader, we repeat the argument: Let $(m_k)_{k\in\mathbb{N}}$ be a sequence of regular values for $f$ such that $m_k \searrow 0$ as $k \to \infty$. Then, $\{f > m_k\} \subset M$ is a regular, open set with outer unit normal $\nu$ in the direction of $-\nabla f$. By Green's formula

$$\begin{aligned}
\int_{\{f>m_k\}} \left(\varphi \Delta_{g_0} f - f \Delta_{g_0} \varphi\right) d\mu_{g_0} &= \int_{\partial\{f>m_k\}} \left(\varphi \langle \nabla f, \nu\rangle_{g_0} - f \langle \nabla \varphi, \nu\rangle_{g_0}\right) d\sigma \\
&\leq -m_k \int_{\partial\{f>m_k\}} \langle \nabla \varphi, \nu\rangle_{g_0} d\sigma \\
&= -m_k \int_{\{f>m_k\}} \Delta_{g_0} \varphi \, d\mu_{g_0} \\
&\leq m_k \int_M |\Delta_{g_0} \varphi| \, d\mu_{g_0}.
\end{aligned}$$

Passing to the limit $k \to \infty$ proves (32) since $(M, g_0)$ is closed.

For any pair $a, b$ of real numbers satisfying $0 < a \leq b$, there holds

$$b^\eta - a^\eta = \frac{b^\eta a - a^{\eta+1}}{a} \leq \frac{b^{\eta+1} - a^{\eta+1}}{a}, \tag{33}$$

$$b^\eta - a^\eta = \frac{b^{\eta+1} - a^\eta b}{b} \leq \frac{b^{\eta+1} - a^{\eta+1}}{b}. \tag{34}$$

Since $w(\cdot, t) > 0$ implies $u(\cdot, t) > \tilde{u}(\cdot, t)$ we may apply (33) with $a = \tilde{u}(x, t)$ and $b = u(x, t)$ and (32) in equation (31) to obtain

$$\frac{1}{\eta+1}\Psi(t) \leq \frac{\sup_M(-\mathrm{R}_{g_0})}{\inf_{M\times[0,T]}\tilde{u}} J(t) + \frac{m-1}{\eta} \int_{\{w(\cdot,t)>0\}} (u^\eta - \tilde{u}^\eta)(\cdot, t) \Delta_{g_0} \varphi \, d\mu_{g_0}.$$

We claim that given $\varepsilon > 0$, there exists $0 \leq \varphi \in \mathrm{C}^\infty_c(M \setminus N)$ satisfying

$$\varphi(p) = \begin{cases} 0, & \text{if } \mathrm{dist}_{g_0}(p, N) < \tfrac{1}{2}\varepsilon, \\ 1, & \text{if } \mathrm{dist}_{g_0}(p, N) > \varepsilon, \end{cases} \tag{35}$$

$$\Delta_{g_0} \varphi \leq C\varepsilon^{-2} \tag{36}$$

with some constant $C$ independent of $\varepsilon$. Indeed, if $r$ denotes the distance from $N$ in $(M, g_0)$ and if $\varepsilon > 0$ is sufficiently small, we may choose $\varphi$ locally depending only on $r$ satisfying (35) as well as $0 \leq \frac{\partial}{\partial r}\varphi \leq Cr\varepsilon^{-2}$ and $\frac{\partial^2}{\partial r^2}\varphi \leq C\varepsilon^{-2}$. Then, (36) follows because Lemma 1.2 implies that $\Delta_{g_0} r \leq \frac{C}{r}$ in the support of $\frac{\partial}{\partial r}\varphi$. Moreover,

$$\int_{\{\Delta_{g_0}\varphi>0\}} \Delta_{g_0} \varphi \, d\mu_{g_0} \leq C\varepsilon^{m-n-2}$$

with some constant $C$ depending on $(M, g_0)$ and $N$. We conclude

$$\frac{1}{\eta+1}\Psi(t) \leq \frac{\sup_M(-\mathrm{R}_{g_0})}{\inf_{M\times[0,T]}\tilde{u}} J(t) + C\varepsilon^{m-n-2} \sup_{M\times[0,T]} u^\eta. \tag{37}$$

Abbreviating

$$\alpha := (\eta + 1)\frac{\sup_M(-\mathrm{R}_{g_0})}{\inf_{M\times[0,T]} \tilde{u}}, \qquad \beta := (\eta + 1)C \sup_{M\times[0,T]} u^\eta$$

and recalling $J(0) = 0$, we obtain

$$J(t) \leq \frac{\beta e^{\alpha t}}{\alpha} \varepsilon^{m-n-2}$$

by applying Lemma A.3 given in the appendix to (37). Since $\varepsilon > 0$ is arbitrary and $m - n - 2 \geq \frac{m}{2} - 1 > 0$, we obtain $J(t) \leq 0$ for every $t \in [0, T]$, or equivalently, $u \leq \tilde{u}$.

Switching the roles of $u$ and $\tilde{u}$ and using estimate (34) instead of (33) shows that the reverse inequality also holds. Therefore, $u = \tilde{u}$. □

**Acknowledgements**  The author is grateful to Professor Alessandro Carlotto and Professor Michael Struwe for inspiring conversations about the incompleteness of Yamabe flows on general punctured manifolds.

# A. Appendix

In the proofs of Theorems 2 and 5, we encounter continuous functions of time which are not quite differentiable. In these cases, we deal with the lim sup of backward difference quotients. The proofs are inspired by Richard Hamilton [8, Lemma 3.1] who obtained similar inequalities for the forward difference quotient.

Lemmata A.1 and A.2 were already proven in [22] and [21, Lemma A.6] but since the proofs are short, we also include them here in this appendix for completeness.

**Lemma A.1.** *Let $T > 0$ and let $f \colon [0, T] \to \mathbb{R}$ be continuous satisfying*

$$\limsup_{\tau \searrow 0} \frac{f(\xi) - f(\xi - \tau)}{\tau} \leq 0$$

*for every $0 < \xi \leq T$. Then $f(t) \leq f(0)$ for every $t \in [0, T]$.*

*Proof.* Let $0 < t \leq T$ and $\varepsilon > 0$ be arbitrary but fixed. By assumption, there exists $\delta > 0$ such that

$$\forall \tau \in [0, \delta[ \colon \quad f(t) - f(t - \tau) \leq \varepsilon \tau. \tag{38}$$

We may assume that $\delta \in \,]0, t]$ is maximal with this property. By continuity of $f$, estimate (38) extends to

$$f(t) - f(t - \delta) \leq \varepsilon \delta. \tag{39}$$

If $t - \delta > 0$, we repeat the argument to find $\delta' > 0$ such that

$$\forall \tau \in [0, \delta'[ \colon \quad f(t - \delta) - f(t - \delta - \tau) \leq \varepsilon \tau. \tag{40}$$

In particular, (39) and (40) can be combined to

$$f(t) - f(t - \delta - \tau) \leq \varepsilon(\delta + \tau)$$

for all $\tau \in [0, \delta'[$ in contradiction to the maximality of $\delta$. Hence, $\delta = t$ and we obtain

$$f(t) - f(0) \leq \varepsilon t.$$

Since $\varepsilon > 0$ is arbitrary, the claim follows. □

**Lemma A.2.** *Let $\eta > 0$ and $Q \geq 0$. Let $v \colon [0, T] \to [0, \infty[$ be a continuous function satisfying*

$$\limsup_{\tau \searrow 0} \frac{v(\xi) - v(\xi - \tau)}{\tau} \leq \eta Q v(\xi)^{1 - \frac{1}{\eta}}$$

*for every $0 < \xi \leq T$ where $v(\xi) \neq 0$. Then $v(t)^{\frac{1}{\eta}} - v(0)^{\frac{1}{\eta}} \leq Qt$ for every $t \in [0, T]$.*

*Proof.* We consider the function $f = v^{\frac{1}{\eta}}$. Let $0 < \xi \leq T$ be arbitrary. If $v(\xi) \neq 0$ then by assumption

$$\eta Q \geq f(\xi)^{1-\eta} \limsup_{\tau \searrow 0} \frac{1}{\tau}\Big(f(\xi)^\eta - f(\xi - \tau)^\eta\Big)$$
$$= \limsup_{\tau \searrow 0} \frac{1}{\tau}\Big(f(\xi) - f(\xi-\tau)\Big)\frac{f(\xi)^\eta - f(\xi-\tau)^\eta}{f(\xi) - f(\xi-\tau)} f(\xi)^{1-\eta}$$
$$= \limsup_{\tau \searrow 0} \frac{\eta}{\tau}\Big(f(\xi) - f(\xi-\tau)\Big). \tag{41}$$

Dividing by $\eta$ yields

$$\limsup_{\tau \searrow 0} \frac{f(\xi) - f(\xi-\tau)}{\tau} \leq Q. \tag{42}$$

If $v(\xi) = 0$ then $f(\xi) = v(\xi)^{\frac{1}{\eta}} = 0$ and we have

$$\limsup_{\tau \searrow 0} \frac{f(\xi) - f(\xi-\tau)}{\tau} \leq 0$$

since $f \geq 0$. We conclude that (42) holds not only where $v(\eta) \neq 0$ but in fact for every $0 < \xi \leq T$ because $Q \geq 0$ is assumed. The claim follows by applying Lemma A.1 to $h(t) = f(t) - Qt$. □

**Lemma A.3.** *Let $a, T > 0$ and $b \geq 0$ and let $J\colon [0,T] \to \mathbb{R}$ be continuous satisfying*

$$\limsup_{\tau \searrow 0} \frac{J(t) - J(t-\tau)}{\tau} \leq aJ(t) + b$$

*for every $0 < t \leq T$. Then, $J(t) \leq J(0)\, e^{at} + (e^{at} - 1)\frac{b}{a}$.*

*Proof.* The map $f\colon [0,T] \to \mathbb{R}$ given by $f(t) = e^{-at}(J(t) + \frac{b}{a})$ is continuous satisfying

$$f(t) - f(t-\tau) = e^{-at}(J(t) + \tfrac{b}{a}) - e^{-at}e^{a\tau}\Big(J(t-\tau) + \tfrac{b}{a}\Big)$$
$$= e^{-at}\Big(J(t) - J(t-\tau)\Big) + e^{-at}\Big(J(t-\tau) + \tfrac{b}{a}\Big)(1 - e^{a\tau})$$

for every $0 \leq \tau \leq t \leq T$. Let $0 < t \leq T$ be arbitrary but fixed. By continuity of $J$,

$$\lim_{\tau \searrow 0}\Big(J(t-\tau) + \tfrac{b}{a}\Big)\frac{1 - e^{a\tau}}{\tau} = \Big(J(t) + \tfrac{b}{a}\Big)(-a) = -aJ(t) - b.$$

Therefore, and by the assumption on the backward difference quotient of $J$, we obtain

$$\limsup_{\tau \searrow 0} \frac{f(t) - f(t-\tau)}{\tau} \leq 0.$$

Lemma A.1 implies $f(t) \leq f(0)$ for every $t \in [0,T]$. Consequently,

$$J(t) + \tfrac{b}{a} \leq e^{at}\Big(J(0) + \tfrac{b}{a}\Big). \qquad \square$$

**Lemma A.4.** *Let $(M, g)$ be a Riemannian manifold of dimension $m \geq 2$ and let $\tilde{g} = ug$ be a conformal metric on $M$. Let $f \colon M \to \mathbb{R}$ be a smooth function. Then, the following relations hold.*

$$|\nabla f|^2_{\tilde{g}} = \tfrac{1}{u}|\nabla f|^2_g \tag{43}$$

$$\Delta_{\tilde{g}} f = \tfrac{1}{u}\Delta_g f + \tfrac{m-2}{2}u^{-2}\langle \nabla u, \nabla f\rangle_g \tag{44}$$

*Proof.* We introduce local coordinates on $M$ and denote by $\tilde{g}_{ij}$ the components of the metric $\tilde{g} = ug$ and by $\tilde{g}^{ij}$ the components of its inverse $\tilde{g}^{-1} = \tfrac{1}{u}g^{-1}$. Using Einstein's summation convention, we have

$$|\nabla f|^2_{\tilde{g}} = \langle \nabla f, \nabla f\rangle_{\tilde{g}} = \tilde{g}_{ij}(\tilde{g}^{\ell i}\partial_\ell f)(\tilde{g}^{jk}\partial_k f) = \tfrac{1}{u}|\nabla f|^2_g,$$

$$\Delta_{\tilde{g}} f = \frac{1}{\sqrt{|\det \tilde{g}|}}\partial_j\left(\sqrt{|\det \tilde{g}|}\,\tilde{g}^{ij}\partial_i f\right)$$

$$= \frac{1}{u^{\frac{m}{2}}\sqrt{|\det g|}}\partial_j\left(u^{\frac{m}{2}-1}\sqrt{|\det g|}\,g^{ij}\partial_i f\right)$$

$$= \tfrac{1}{u}\Delta_g f + \tfrac{m-2}{2}u^{-2}\langle \nabla u, \nabla f\rangle_g. \qquad \square$$